\documentclass[11pt,letterpaper]{article}
\RequirePackage[OT1]{fontenc}
\RequirePackage[noinfoline,amsthm,amsmath]{imsart}
\usepackage[authoryear]{natbib}
\allowdisplaybreaks

\startlocaldefs
\numberwithin{equation}{section}
\theoremstyle{plain}
\newtheorem{thm}{Theorem}[section]
\newtheorem{lemma}{Lemma}[section]
\newtheorem{corollary}{Corollary}[section]
\newtheorem{example}{Example}[section]
\endlocaldefs
\begin{document}

\begin{frontmatter}
\title{
Admissibility and minimaxity of generalized Bayes estimators for
spherically symmetric family}
\runtitle{Admissibility for spherical family}

\begin{aug}
\author{\fnms{Yuzo} \snm{Maruyama}
 \thanksref{t1}
\ead[label=e1]
{maruyama@csis.u-tokyo.ac.jp}
}
and
\author{\fnms{Akimichi} \snm{Takemura}
\ead[label=e2]
{takemura@stat.t.u-tokyo.ac.jp}
}

\thankstext{t1}{Corresponding author}
\runauthor{Y.Maruyama and A.Takemura}

\affiliation{The University of Tokyo}

\address{
Center for Spatial Information Science, \\
The University of Tokyo\\
5-1-5 Kashiwanoha, Kashiwa-shi, Chiba 
277-8568, Japan\\
E-mail: \printead[maruyama@csis.u-tokyo.ac.jp]{e1}}

\address{
Graduate School of Information Science and Technology\\
The University of Tokyo\\
7-3-1 Hongo, Bunkyo-ku, Tokyo 
113-0033, Japan\\
E-mail: \printead[takemura@stat.t.u-tokyo.ac.jp]{e2}}

\end{aug}

\begin{abstract}
  We give a sufficient condition for admissibility of generalized Bayes
  estimators of the location vector of spherically symmetric distribution
  under squared error loss.  Compared to the known results for the
  multivariate normal case, our sufficient condition is very tight and is
  close to being a necessary condition.  In particular we establish the
  admissibility of generalized Bayes estimators with respect to the
  harmonic prior and priors with slightly heavier tail than the harmonic
  prior.  We use the theory of regularly varying functions to construct a
  sequence of smooth proper priors approaching an improper prior fast
  enough for establishing the admissibility. We also discuss conditions of
  minimaxity of the generalized Bayes estimator with respect to the
  harmonic prior.
\end{abstract}



\end{frontmatter}

\section{Introduction}
\label{sec:intro}
We consider estimation of the $p$-dimensional location parameter of a
spherically symmetric distribution. Specifically, let $ X =
(X_1,\ldots,X_p)' $ have a density function $ f(\|x-\theta \|) $ and
consider estimation of $ \theta $ with a general quadratic loss
function $ L_Q(\theta,\delta)= ( \delta - \theta)'Q( \delta -
\theta)=\|\delta - \theta \|^2_Q $ for a positive definite matrix $Q$.
The usual minimax estimator $ X $, which is generalized Bayes, is
inadmissible for $ p \geq 3 $ as shown in \cite{Stein1956} in the normal
case, and in \cite{Brown1966} under more general situation, respectively.
In the decision-theoretic point of view, we are interested in
proposing admissible estimators dominating $ X $, that is, minimax
admissible estimators.  Note that the dominance over $X$ means
minimaxity in our setting because $ X$ is minimax with a constant
risk.  Note also that our results hold for the location vector of
elliptically-contoured distributions, because we are considering a
general quadratic loss function with arbitrary positive definite $Q$.

In the normal case, there already exists a broad class of admissible
minimax estimators. \cite{Baranchik1970} gave a sufficient condition for
minimaxity of a shrinkage estimator of the form
\begin{equation} 
\label{orthogonal}
 \delta_{\phi}(X)=(1-\phi(\|X\|)\|X\|^{-2})X.
\end{equation}
\cite{Strawderman1971} found a subclass of proper Bayes estimators
of the form \eqref{orthogonal} 
satisfying the sufficient conditions for minimaxity.
\cite{Brown1971} gave a very powerful sufficient condition for
admissibility of generalized Bayes estimators.
Using Brown's (1971) condition, 
\cite{Berger1976b}, \cite{Fourdrinier-etal1998} and
\cite{Maruyama1998, Maruyama2004}
enlarged a class of admissible minimax estimators
which are of the Strawderman type and generalized Bayes. 
For a subclass of scale mixtures of multivariate normal distributions
which includes multivariate-$t$ distribution, 
some proper Bayes minimax estimators were proposed by \cite{Maruyama2003}
by using Strawderman's (1971) techniques.

However, for general spherically symmetric distributions, no minimax
admissible estimators of the location vectors have been derived,
although for the minimaxity, various Baranchik-type sufficient
conditions of the estimator \eqref{orthogonal} were given by
\cite{Berger1975}, 
\cite{Brand-Straw1978, Brand-Straw1991b} and \cite{Bock1985}.  The
main reason is the lack of a standard class of generalized or proper
Bayes estimators of the form \eqref{orthogonal} like the Strawderman
type in the normal case, which allows an easy check of the minimaxity
condition.  Furthermore no sufficient condition for admissibility of
generalized Bayes estimators has been derived.  In 
\cite{Maruyama-Takemura2005}, we 
provided satisfactory solutions to these problems.
In this paper we will weaken 
some regularity conditions assumed in \cite{Maruyama-Takemura2005}.

In Section 2, we give preliminary results including the properties of
regularly varying functions and asymptotic behaviors of expected
values when $\|\theta\|$ is sufficiently large.  The former is useful
for constructing a very convenient sequence of proper densities
approaching an improper
density $g(\theta)$, which is required in applying the method of 
\cite{Blyth1951}.  

In Section 3, we will present a powerful sufficient
condition for admissibility of generalized Bayes estimator and in
particular show that the generalized Bayes estimators with respect to
the harmonic prior $g(\theta)=\| \theta \|^{2-p}$ and with respect to 
a prior with a slightly heavier tail 
\begin{align} 
\label{eqn:2-p+log}
g(\theta)= \| \theta \|^{2-p}\log(\|\theta\|+c), \quad c>1, 
\end{align}
are admissible under mild regularity conditions on $f$.  

In Section 4, we show that the generalized Bayes estimator with respect to
the harmonic prior is written as
\begin{equation} \label{harmonic-Bayes}
\left(1-\frac{\int_{0}^{1}t^{p-1}F(t \|X\|)dt}
{\int_{0}^{1}t^{p-3}F(t \|X\|)dt} \right)X,
\end{equation} 
where $ F(u)=\int_{u}^{\infty}sf(s)ds $.  This form is simple enough
to check various sufficient conditions for minimaxity and 
we demonstrate that \eqref{harmonic-Bayes} is minimax for some $f$.
We believe that \eqref{harmonic-Bayes} is 
minimax for a broad subclass of spherically
symmetric distributions.
Notice that the generalized Bayes estimators with respect to
priors except $\|\theta\|^{2-p}$ do not have such simple forms 
as far as we know.

Our proof of admissibility is mainly based on the techniques of
\cite{Brown-Hwang1983}.  \cite{Brown-Hwang1983} considered the problem of
estimating the natural mean vector of an exponential family under a
quadratic loss function.  Note that the intersection of their setting and
our setting is the multivariate normal case.  Their sufficient
condition for admissibility in the normal case does not however permit
$g(\theta)$ to diverge to infinity around the origin like $ \| \theta
\|^{2-p}$,  while it permits $g(\theta) \leq \| \theta\|^a$ with $a\leq 2-p$
for sufficiently large $\|\theta \|$.  Prior to \cite{Brown-Hwang1983},
\cite{Brown1971} considered the estimation in the multivariate normal case
and gave a powerful sufficient condition for minimaxity which are satisfied
by the harmonic prior and \eqref{eqn:2-p+log}, but his proof was based on many advanced
mathematics. Our mathematical tool is much more familiar to the readers.
\cite{Brown1971} also gave sufficient condition for inadmissibility.  By
using it we see that the generalized Bayes estimator with respect to
$\| \theta \|^{2-p}\log^2 (\|\theta\|+2)$ is inadmissible. Hence our
sufficient condition for admissibility should be very tight and close to
being a necessary condition.  

\cite{Brown1979} considered a more general problem than ours: estimation of
$ \theta $ for a general density $p(x-\theta)$ and a general loss function
$ W(\delta-\theta)$.  He conjectured that the prior $g(\theta) \sim
\|\theta\|^{a}$ with $a\leq 2-p$ for sufficiently large $ \|\theta\|$
leads to admissibility, regardless of the density $p$ and the loss $W$.
Hence our results support \cite{Brown1979}'s conjecture for the case of
elliptically-contoured family and a general quadratic loss function.

Finally we notice that the most important key for our proof for
admissibility is the construction of a very convenient sequence
$h_i(\theta)$ for approximating $g(\theta)$ by 
$g(\theta)h^\gamma_i(\theta)$ for $0<\gamma \le 2$.
\cite{Brown-Hwang1983} used $ \gamma=2$ and
\begin{align*}
 h_i(\theta)=\begin{cases}
	      1 & \|\theta\| \leq 1 \\
1-\log \|\theta\|/\log i & 1 \leq \| \theta\| \leq i \\
0 & \|\theta\|>i .
	     \end{cases} 
\end{align*} 
This $ h_i(\theta) $ is not differentiable at $ \|\theta\|=1$ and
truncated at $\|\theta\|=i$, which makes handling and extension difficult
for our purposes.  Our $h_i(\theta)$ given in Section 2 is smoother and
not truncated. 
Furthermore our $\gamma$ is flexible whereas $\gamma=2$ in 
\cite{Brown-Hwang1983}. By such a flexible $\gamma$, 
we can adjust the rate of convergence of
$ g(\theta)h^\gamma_i(\theta) $ so that it is just enough to be proper.
We will see that choosing as small $\gamma$ as possible is
important in the main theorem, Theorem 3.1.
We believe that our smooth function $h_i(\theta)$ and an idea of 
flexible $\gamma$ are very useful for showing
admissibility of generalized Bayes estimators in various problems.

\section{Preliminaries}
\label{sec:preliminaries}
In this section we prepare a sequence of proper densities using the
theory of regularly varying functions and give some results on
asymptotic behaviors of expected values when the location parameter
diverges to infinity.  For the theory of regularly varying and slowly
varying functions the readers are referred to \cite{geluk-dehaan1987}
and \cite{BGT1987}.

\subsection{Regularly varying functions}
\label{sec:rv}
A Lebesgue measurable function $ f: R^+ \to R $ which is eventually
positive is called regularly varying if for some $\alpha \in R$
\begin{align} 
\label{def-RV}
 \lim_{x \to \infty}\frac{f(tx)}{f(x)}=t^\alpha , \quad  \forall t>0. 
\end{align}
We sometimes use the notation $ f \in \mbox{RV}_\alpha$.
The number $ \alpha $ in the above is called the index of regular
variation. A function satisfying \eqref{def-RV} with $\alpha=0$ is
called slowly varying.

Let $ \beta: R^+ \to R^+ $ be positive, continuously
differentiable, monotone decreasing, integrable (i.e.\ $ \int_0^\infty
\beta(r)dr<\infty$) and regular varying with index $-1$.  A typical $
\beta(\eta)$ is
\begin{align} 
\label{beta1}
\frac{1}{\eta+c}\frac{1}{\{\mathrm{Log}_n(\eta+c)\}^2}
\prod_{i=0}^{n-1}\frac{1}{\mathrm{Log}_i(\eta+c)},
\end{align}
where $n$ is a positive integer, $\mathrm{Log}_0(\eta+c) \equiv 1 $,
\begin{align*}
\mathrm{Log}_i(\eta+c)=\underbrace{\log \log \cdots \log}_{i}(\eta+c),
\quad i \ge 1,
\end{align*}
and $c$ is chosen such that $  \mathrm{Log}_n(c)>0 $.
Note that for
\eqref{beta1}
\[
 \int_\eta^\infty \beta(r)dr= \frac{1}{\mbox{Log}_n(\eta+c)}.
\]

The following results for $\beta(\eta)$ satisfying the above
assumptions are known from the theory of regularly varying functions.
\begin{lemma} \label{karamata}
\begin{enumerate}
 \item $ \int_{\eta}^\infty \beta(r)dr \in \mbox{\rm RV}_0$ and 
$ \beta'(r) \in \mbox{\rm RV}_{-2}$.
\item $ \lim_{\eta \to \infty} \eta \beta(\eta)/\int_{\eta}^\infty
       \beta(r)dr=0$, 
$ \lim_{\eta \to \infty} \eta \beta'(\eta)/\beta(\eta)=-1$.
\end{enumerate} 
\end{lemma}

We now define  functions $H_i(\eta)$, $i=1,2,\ldots$, based on $
\beta(r)$ by 
\begin{align} 
\label{H}
H_i(\eta)=
\frac{\int_\eta^\infty e^{ (\eta-r)/i}\beta(r)dr} 
{\int_\eta^\infty \beta(r)dr}.
\end{align}
These functions are very useful for constructing a sequence of proper
prior densities approaching the target improper density in the next
section.  The properties of $H_i$ are given in the following theorem.

\begin{thm} 
\label{thm-beta-H}
 \begin{enumerate}
\item  $0 \le H_1(\eta) \le H_2(\eta) \le \cdots \le 1$.  For 
any fixed $ \eta$, 
$ \lim_{i \to \infty}H_i(\eta)=1 $.
\item For any fixed $i$, $ \lim_{\eta \to \infty}
\int_\eta^\infty \beta(r)dr \beta(\eta)^{-1}
H_i(\eta)=i$ and hence $ H_i(\eta) \in \mbox{\rm RV}_{-1}$.
\item For any fixed $\eta$, $ \lim_{i \to
	\infty}H'_i(\eta)=0 $.
\item $  | H'_i(\eta) |  < 2 \beta(\eta)/\int_\eta^\infty
\beta(r)dr$ for all $\eta > 0$.
\item For any $ \epsilon>0$, there exists $ \eta_0$ such that
$ -1 -\epsilon < \eta H'_i(\eta)/H_i(\eta) \leq 0 $ for 
all $\eta \geq \eta_0$ and for all $i$.
\end{enumerate}
\end{thm}

\begin{proof}
It is obvious that $0 \le H_1(\eta)  \le 1$  and $H_i(\eta)$ is
increasing in $i$. For fixed  $  \eta $, 
$H_i(\eta) \uparrow 1$ by the monotone convergence theorem. 

By integration by parts, the numerator of $H_i(\eta)$ is written as
\begin{equation}
\label{eq:numerator-integ-parts}
\int_\eta^\infty e^{ (\eta-r)/i}\beta(r)dr = 
i \beta(\eta) + i \int_\eta^\infty e^{ (\eta-r)/i}\beta'(r)dr.
\end{equation}
Therefore 
\begin{align}
 H_i(\eta)=i \frac{\beta(\eta)}{\int_\eta^\infty \beta(r)dr} 
+i
\frac{\int_\eta^\infty e^{ (\eta-r)/i}\beta'(r)dr}
{\int_\eta^\infty \beta(r)dr}.
\label{H+}
\end{align}
\eqref{H+} divided by \eqref{H} is
\begin{align*}
 1 = 
i \frac{\beta(\eta)}{H_i(\eta)\int_\eta^\infty \beta(r)dr} 
+i
\frac{\int_\eta^\infty e^{ -r/i}\beta'(r)dr}
{\int_\eta^\infty e^{ -r/i}\beta(r)dr} .
\end{align*}
For fixed $i$, the second term of the above equation converges to $0$
as $ \eta \rightarrow \infty $ by the L'Hospital theorem.  $H_i(\eta)
\in \mbox{RV}_{-1}$ because $ \int_\eta^\infty \beta(r)dr \in
\mbox{RV}_0$ and $ \beta(\eta) \in \mbox{RV}_{-1}$.

Using (\ref{eq:numerator-integ-parts}) again, differentiation of
the numerator of $H_i(\eta)$ gives
\[
\left( \int_\eta^\infty e^{ (\eta-r)/i}\beta(r)dr\right)'=
\frac{1}{i}\int_\eta^\infty e^{ (\eta-r)/i}\beta(r)dr- \beta(\eta)
= \int_\eta^\infty e^{ (\eta-r)/i}\beta'(r)dr .
\]
Therefore
\begin{align}
 H'_i(\eta) 
=
\frac{\beta(\eta)\int_\eta^\infty e^{ (\eta-r)/i}\beta(r)dr}
{(\int_\eta^\infty \beta(r)dr)^2}
-\frac{\int_\eta^\infty e^{ (\eta-r)/i}\{-\beta'(r)\}dr}
{\int_\eta^\infty \beta(r)dr}.
 \label{bibun-H}  
\end{align}
Note that $-\beta'(r)\ge 0$ by our assumption.  Each term of the right
hand side of \eqref{bibun-H} is nondecreasing in $i$ and hence by the
monotone convergence theorem
\begin{align*}
 \lim_{i \to \infty}H'_i(\eta)
= 
\frac{\beta(\eta)\int_\eta^\infty \beta(r)dr}
{(\int_\eta^\infty \beta(r)dr)^2}
-\frac{\int_\eta^\infty \{-\beta'(r)\}dr}
{\int_\eta^\infty \beta(r)dr}
 =0 .
\end{align*}
Furthermore we have
\begin{eqnarray*}
 \left| H'_i(\eta)\right| 
& < &
\left|
\frac{\beta(\eta)\int_\eta^\infty e^{ (\eta-r)/i}\beta(r)dr}
{(\int_\eta^\infty \beta(r)dr)^2}\right|
+\left|\frac{\int_\eta^\infty e^{ (\eta-r)/i}\{-\beta'(r)\}dr}
{\int_\eta^\infty \beta(r)dr}\right|\\
& <& 2\beta(\eta)/\int_\eta^\infty \beta(r)dr.
\end{eqnarray*}

Dividing \eqref{bibun-H} by \eqref{H}, we have
\begin{align}
 \eta \frac{H'_i(\eta)}{H_i(\eta)}
&=\eta \left(\frac{\beta(\eta)}{\int_\eta^\infty
\beta(r)dr}+\frac{\int_\eta^\infty e^{ -r/i}\beta'(r)dr}
{\int_\eta^\infty e^{-r/i}\beta(r)dr}\right) \label{bibun-HH} \\
&>  \frac{\eta\beta(\eta)}{\int_\eta^\infty
\beta(r)dr}- \frac{\int_\eta^\infty e^{ -r/i}\{-r \beta'(r)/\beta(r)\}
\beta(r)dr}
{\int_\eta^\infty e^{-r/i}\beta(r)dr} \nonumber \\ 
&>  \frac{\eta\beta(\eta)}{\int_\eta^\infty
\beta(r)dr}- \sup_{r>\eta}\frac{-r \beta'(r)}{\beta(r)}. \nonumber
\end{align}
By 2 of Lemma \ref{karamata} the right hand side converges to $-1$.
This implies that for any $ \epsilon>0$ there exists $ \eta_0$
such that $ \eta H'_i(\eta)/H_i(\eta)>-1-\epsilon$ for all  
$\eta \geq \eta_0$ and for all $i$.
Finally 
we will prove that $H'_i(\eta) \leq
0$ for sufficiently large $\eta$ independent of $i$. 
 By 2 of Lemma \ref{karamata}, 
\begin{align*}
 \frac{\eta \int_\eta^\infty \beta(r)dr}{\beta(\eta)}
\left(\frac{\beta(\eta)}{\int_\eta^\infty \beta(r)dr}\right)'
=\eta \frac{\beta'(\eta)}{\beta(\eta)}+ \eta \frac{\beta(\eta)} 
{ \int_\eta^\infty  \beta(r)dr} \to -1
\end{align*}
and hence
$ \beta(\eta)/\int_\eta^\infty
\beta(r)dr$ is eventually nonincreasing.
Hence by redefining $\eta_0$ if necessary, we can assume
that $ \beta(\eta)/\int_\eta^\infty
 \beta(r)dr$ is monotone nonincreasing for $\eta
\ge \eta_0$. 
By integration by parts on the numerators of each term in
\eqref{bibun-H}, we have
\begin{align*}
& \int_\eta^\infty e^{ -r/i}\beta(r)dr= e^{-\eta/i}\int_\eta^\infty 
\beta(r)dr -i^{-1} \int_\eta^\infty e^{ -r/i} \left\{\int_r^\infty 
\beta(s)ds\right\}dr \\
& \int_\eta^\infty e^{ -r/i}\beta'(r)dr= -e^{-\eta/i}\beta(\eta)
+i^{-1} \int_\eta^\infty e^{ -r/i} \beta(r)dr 
\end{align*}
and hence
\begin{align*}
\left\{  \frac{i \int_\eta^\infty  \beta(r)dr }{\int_\eta^\infty e^{
 -r/i} \left\{\int_r^\infty  
\beta(s)ds\right\}dr} \right\}
H'_i(\eta)
= -\frac{\beta(\eta)}{\int_\eta^\infty \beta(r)dr}+
\frac{\int_\eta^\infty e^{ -r/i} \beta(r)dr}{\int_\eta^\infty e^{
 -r/i} \left\{\int_r^\infty  \beta(s)ds\right\}dr},
\end{align*}
which is nonpositive for $ \eta \ge \eta_0$.
%
\end{proof}

\subsection{Asymptotic behavior of expectations}
\label{sec:asymptotic-behavior}
In the next section, we need evaluation of an asymptotic
behavior of expectation
\begin{align*}
 E_x[\rho(\theta)]=\int_{R^p} \rho(\theta)f(\|\theta - x \|)d\theta
\end{align*}
for sufficiently large $\|x\|$, where a random vector $ \theta$ has
the density function $ f(\|\theta - x\|)$.  This is the expected value
with respect to the posterior distribution. Interchanging the roles of
$x$ and $\theta$, in this subsection, 
we consider the asymptotic behavior of expectation
\begin{align*}
 E_\theta[\rho(X)]=\int_{R^p} \rho(x)f(\| x - \theta \|)dx
\end{align*}
for sufficiently large $\|\theta\|$, where a random vector $ X$ has the
density function $ f(\|x - \theta \|)$.  We believe that this does not
confuse the readers.  

We discuss some notations used in the following. In addition to the
Euclidean norm $\|x\|^2=x_1^2 + \cdots + x_p^2$, we consider the norm
$\|x\|_d^2=d_1^2 x_1^2 + \cdots + d_p^2 x_p^2$.  For convenience we assume,
without loss of generality, that $ d_1 \ge \cdots \ge d_p \ge 1.  $ Under
this assumption
\begin{equation}
\label{eq:norm-d}
 \| x \| \le \| x \|_d \le d_1 \| x \|.
\end{equation}
By introducing this norm our results hold for
elliptically-contoured distributions.  
The gradient of $\rho(x)$ is denoted by
\[
\nabla \rho(x)=(\frac{\partial}{\partial x_1} \rho(x),\ldots,\frac{\partial}{\partial
x_p}\rho(x))' .
\]
We also write $\nabla_j \rho(x)=(\partial/\partial x_j) \rho(x)$.  Finally 
we write $c_p=2\pi^{p/2}/\Gamma(p/2)$.

Now we make the following regularity conditions on the density $f$ and
the function $\rho $. 
\begin{description}
\item[F1] There exist $r_0>0$, $L>0$, and $s>1$, such that
  $r^{p+s}f(r) \leq L$ for all $r\geq r_0$.
\item[B1] 
  $\rho(x)$ is written as $ \rho(x)=\varrho(\| x \|_d)$, where 
$\varrho(r)$ is continuously differentiable in $r>0$.
\item[B2]  There exists $ r_1 \ge 1$ and $ t_1 \leq t_2$
such that
$ \varrho(r)>0 $ and $ t_1 \leq r\varrho'(r)/\varrho(r) \leq t_2$ for
	   all $r \geq r_1$.
\end{description}
Assumption {\bfseries B2} is, for instance, satisfied by 
\begin{align*}
 \varrho(r)=\exp\left( \int_0^r \frac{-p+\alpha\cos t}{t+1} dt \right), 
\mbox{ for }\alpha>0
\end{align*}
where we easily see 
$ r\varrho'(r)/\varrho(r)=\{r/(r+1)\}(-p+\alpha \cos r)$ and hence
$t_1=-p-\alpha$ and $t_2=-p+\alpha$ in {\bfseries B2}.
Since $ \lim_{r \to \infty} r\varrho'(r)/\varrho(r)$ exists for 
regularly varying $\varrho$, 
we deal with a broader class of functions
than the class of regularly varying functions.
We will discuss more in Section 3.
Note that
\[
\nabla \rho(x)=\frac{\varrho'(\| x \|_d)}{\| x \|_d} 
(d_1^2 x_1, \ldots, d_p^2 x_p)'
\]
and
\[
\| \nabla \rho(x)\| = \frac{|\varrho'(\| x \|_d)|}{\| x \|_d} (d_1^4 x_1^2 +
\cdots + d_p^4 x_p^2)^{1/2}\le d_1 |\varrho'(\| x \|_d)|.
\]

The following lemma is useful. The proof based on the integration of
$(\log \varrho(r))'=\varrho'(r)/\varrho(r)$ is easy and omitted.
\begin{lemma}
\label{lem:b3}
Under the assumption {\bfseries B2}
\begin{align*}
 (z/y)^{t_1}\leq \varrho(z)/\varrho(y) \leq (z/y)^{t_2}
\end{align*}
for any $ z>y \geq r_1$. Moreover 
\[
\limsup_{y \to \infty}\sup_{\alpha y \leq z \leq \beta y} 
\varrho(z)/\varrho(y) \le 
\max(\alpha^{t_1},\beta^{t_2})\]
for any $ 0<\alpha<1<\beta $.
\end{lemma}

We now state the following theorem concerning the asymptotic behavior
of $E_\theta[\rho(X)]$ for large $\|\theta\|_d$.

\begin{thm} 
\label{thm-zenkin}
Assume {\bf F1}, {\bf B1} and {\bf B2}.  
For $a=0$ or $1$, and $j=1,\dots,p$, 
if $ s> \max(1, -t_1-a-p, t_2+a)$ and
$ \int_0^1 r^{p+a-1}| \varrho(r)|dr <\infty$,
then there exists $ \epsilon>0$ (say 
$\epsilon=  \min(1, s+t_1+a+p)/4$)
 such that
\begin{equation}
\label{eq:thm22}
\|\theta\|_d^{\epsilon-a}
\left|E_\theta[X_j^a \rho(X)] - \theta_j^a \rho(\theta)
 \right| 
< C \rho(\theta)
\end{equation}
for  $\|\theta\|_d \geq 2d_1\max(r_0,r_1)$.
Moreover $C$ depends on $ \rho$ (or $ \varrho$) only through
$r_1$, $t_1$, $t_2$ and
$ \{\varrho(r_1)\}^{-1}\int_0^{r_1}r^{p+a-1}|\varrho(r)|dr$.
\end{thm}
For simplicity, in the rest of the paper, we will write 
$ E_\theta[X_j^a \rho(X)] \approx \theta_j^a \rho(\theta)$
if, as in Theorem \ref{thm-zenkin}, there exists $\epsilon>0$ such that
\eqref{eq:thm22} is satisfied for sufficiently large $\|\theta\|_d$.
\begin{proof}
Note that 
\[E_\theta[X_j^a \rho(X)] - \theta_j^a \rho(\theta)=
E_\theta[X_j^a (\rho(X)-\rho(\theta))].
\]
 Fix $0 < \nu < 1$ (set $\nu=1/2$ finally).
Define
\begin{eqnarray*}
V_{\nu}&=&\{ x: \|x-\theta \| \leq \nu\| \theta \| \} \\
V'_{\nu}&=&\{ x: (1-\nu)\|\theta\| \leq \|x\| \leq (1+\nu) \|\theta \| \}\\
V''_{\nu}&=&\{ x: d_1^{-1}(1-\nu)\|\theta\|_d \leq \|x\|_d \leq
d_1 (1+\nu) \|\theta \|_d \} .
\end{eqnarray*}
By (\ref{eq:norm-d}) $ V_\nu \subset V'_\nu \subset V''_\nu$.
Then 
\begin{align}
& 
\|\theta\|_d^{\epsilon-a}
\left|E[X_j^a (\rho(X)-\rho(\theta))] \right| \nonumber \\
&  \leq \|\theta\|_d^{\epsilon-a}
\left( \int_{V_\nu} + \int_{V_\nu^C} \right) \|x\|^a 
\, |\rho(x)-\rho(\theta) | \,
f(\|x -\theta\|)dx 
\nonumber \\ 
&  
\leq \|\theta\|_d^{\epsilon}
\sup_{x \in V''_\nu}
\left(\frac{\|x\|_d}{\|\theta\|_d}\right)^a
\int_{V_\nu} | \rho(x) - \rho(\theta)| \,
f(\|x-\theta\|)dx
\nonumber \\
&\quad + \|\theta\|_d^{\epsilon-a} \rho(\theta) \int_{V_\nu^C}
\|x\|^af(\|x -\theta\|)dx 
+ \|\theta\|_d^{\epsilon-a} 
\int_{V_\nu^C} \|x\|^a | \rho(x)|  
f(\|x -\theta\|)dx 
\nonumber \\ 
& =  I_1 + I_2 + I_3.  \quad ({\rm say})
\label{eq:intergral-3-parts}
\end{align}

Consider the first integral $I_1$.  
$ \sup_{x \in V''_\nu}
\left(\|x\|_d/\|\theta\|_d \right)^a  \leq d_1^{a}(1+\nu)^{a}$.
If  $s>1$, then
$m_1=\int_{R^p} \|x-\theta\|f(\|x - \theta\|)dx$ is finite.
Therefore for $ \|\theta\|_d \geq d_1(1-\nu)^{-1}r_1$  we have
\begin{align*}
  & \|\theta\|_d^{\epsilon} \int_{V_\nu} \,|
  \rho(x)-\rho(\theta)| \,
  f(\|x - \theta\|)dx \\
  &\qquad = \|\theta\|_d^{\epsilon} \int_{V_\nu} |(x-\theta)'\nabla
  \rho(x^*)|f(\|x - \theta\|)dx , \
  x^* \in V_\nu \\
  &\qquad \leq m_1 \|\theta\|_d^{\epsilon}
  \sup_{x \in V_\nu} \| \nabla \rho(x) \| \\
  &\qquad \leq m_1 d_1 \|\theta\|_d^{\epsilon} \sup_{x \in
    V_\nu}
  |\varrho'(\|x\|_d)| \\
  &\qquad \leq m_1d_1 \|\theta\|_d^{\epsilon-1} \sup_{x \in
    V''_\nu}\frac{\| \theta \|_d}{\|x\|_d} \sup_{x \in
    V''_\nu}\frac{\varrho(\|x\|_d)}{\varrho(\|\theta\|_d)} \sup_{x \in
    V''_\nu} \frac{\|x\|_d | \varrho'(\|x\|_d) |}{\varrho(\|x\|_d)}
  \times \rho(\theta)
  \\
  & \qquad \le 
  \frac{m_1 d_1^2}{1-\nu}
 \max\left( \left(\frac{1-\nu}{d_1}\right)^{t_1}, 
\left\{d_1(1+\nu)\right\}^{t_2} \right) 
\max(|t_1|, |t_2|)
\times \rho(\theta)
\end{align*}
for $ 0<\epsilon<1$.
Therefore we have $I_1 \le C_{1}\rho(\theta)$, where
$C_{1}= d_1^a(1+\nu)^{a} 
 m_1 \{d_1^2/(1-\nu)\}
\max( \{(1-\nu)/d_1\}^{t_1},
\{d_1(1+\nu)\}^{t_2}) \max(|t_1|, |t_2|)$.

Now we consider the integral outside of $V_\nu$.  We only consider 
$\|\theta\|_d \ge \max(d_1 \nu^{-1} r_0,r_1)$.  
Then  for $x\in V_\nu^C$
\[
\| x- \theta \| \ge \nu \| \theta\| \ge \nu d_1^{-1} \| \theta\|_d \ge
r_0.
\]
Therefore we have, for $0\le \alpha < s$
\begin{align}
  \int_{V_\nu^C}\|x\|^\alpha
  f(\|x - \theta \|)dx 
&\leq  \int_{V_\nu^C}\{\|x-\theta \|+\|\theta\|\}^\alpha
  f(\|x - \theta \|)dx \notag \\
& \leq (1+1/\nu)^{\alpha}\int_{V_\nu^C}\|x-\theta \|^\alpha
  f(\|x - \theta \|)dx \notag \\
& \leq 
(1+1/\nu)^{\alpha}c_p L  \int_{\nu \|\theta\|}^\infty
r^{-s+\alpha-1} dr \notag \\
&=
(1+1/\nu)^{\alpha}c_p L  \frac{(\nu \|\theta\|)^{-s+\alpha}}{s-\alpha}
\notag\\
&\le 
C_2(\alpha)\|\theta\|_d^{\alpha- s},
\label{soto-sekibun}
\end{align}
where $ C_2(\alpha)=(1+1/\nu)^{\alpha}c_p L (d_1/\nu)^{s-\alpha}
 (s-\alpha)^{-1}$.
Hence for the second term $I_2$, 
if $s>1$ and $0<\epsilon<1$, 
then  $I_2 \le C_{2}(a) \rho(\theta)$
 for $a=0,1$.

We have seen that $I_1$ and $I_2$ are bounded from above assuming only $s>1$.

The third term $I_3$  of (\ref{eq:intergral-3-parts}) is more
problematic.  
Write 
\begin{eqnarray*}
I_3 &=&
\|\theta\|_d^{\epsilon-a}  \int_{V_\nu^C}
\|x\|^a|\rho(x)|f(\|x -\theta\|)dx \\
&\le &  \|\theta\|_d^{\epsilon-a}
 \left( 
 \int_{V_\nu^C \cap \{\|x\|_d < r_1 \}} 
 + \int_{V_\nu^C \cap  \{r_1 \leq \|x\|_d \leq \|\theta\|_d \}}  \right.
\\ 
  && \qquad \qquad  \qquad\qquad 
  \left. 
  +\int_{V_\nu^C \cap \{\|x\|_d > \|\theta \|_d \} } \right) \|x\|^a
  |\rho(x)|  f(\|x -\theta\|)dx \\
&=& I_{31} + I_{32} + I_{33}. \qquad ({\rm say})
\end{eqnarray*}
We take care of $I_{33}$ first.  
Since $\varrho(r)r^{-t_2}$ is monotone nonincreasing for $ r \ge r_1$, 
$\rho(x)\|x\|_d^{-t_2} \le \rho(\theta)\| \theta \|_d^{-t_2}$ 
for $\|x\|_d > \|\theta \|_d(\ge r_1)$.
Therefore 
we have, 
\begin{align*}
I_{33} & \le 
\|\theta\|_d^{\epsilon-a-t_2}  \rho(\theta) d_1^{\max(t_2,0)}
\int_{V_\nu^C \cap \{\|x\|_d > \|\theta \|_d \} } \|x\|^{a+t_2} 
f(\|x -\theta\|)dx .
\end{align*}
If $ 0 \leq a+t_2<s$, as in \eqref{soto-sekibun}
\begin{align*}
\int_{V_\nu^C \cap \{\|x\|_d > \|\theta \|_d \} } \|x\|^{a+t_2} 
f(\|x -\theta\|)dx &\leq \int_{V_\nu^C} \|x\|^{a+t_2} 
f(\|x -\theta\|)dx \\ & \leq C_2(a+t_2) \|\theta \|_d^{a+t_2-s}
 \end{align*} 
and if $ a+t_2<0$,
\begin{align*}
\int_{V_\nu^C \cap \{\|x\|_d > \|\theta \|_d \} } \|x\|^{a+t_2} 
f(\|x -\theta\|)dx &\leq d_1^{-t_2-a} \|\theta\|^{a+t_2}_d 
\int_{V_\nu^C} 
f(\|x -\theta\|)dx \\ & \leq d_1^{-t_2-a} C_2(0) \|\theta \|_d^{a+t_2-s}.
 \end{align*}
%
Hence $I_{33}\leq C_{33}\rho(\theta)$ where $ C_{33}= d_1^{\max(t_2,0)}
\max( C_2(a+t_2), C_2(0)d_1^{-a-t_2})$.

Next we consider $I_{31}$.
For $\|\theta\|_d \ge \max(d_1 \nu^{-1}
r_0,r_1)$ and $x\in V_\nu^C$
\begin{align}
f(\| x- \theta \|) \le  L \| x- \theta \|^{-p-s} 
\le L(\nu \|\theta\|)^{-p-s} 
\leq L(d_1/\nu\| \theta\|_d)^{p+s}. \label{f-inequality}
\end{align}
Therefore
\[
I_{31}\le 
\|\theta\|_d^{\epsilon-a} 
Ld_1^{p+s}\nu^{-p-s}\| \theta\|_d^{-p-s}
\int_{\|x\|_d \leq r_1}  \|x\|_d^a \, \rho(x) dx.
\]
Note that by simple change of variables we have
\[
\frac{\partial}{\partial r} \int_{ \|x\|_d \le r} dx = c_{p,d} r^{p-1},
\qquad c_{p,d}=c_p \prod_{i=1}^p d_i^{-1} .
\]
Then
\[
\int_{\|x\|_d \leq r_1}  \|x\|_d^a \, |\rho(x)| dx 
= c_{p,d} \int_0^{r_1} r^{p+a-1}
|\varrho(r)| dr.
\]
Therefore 
\[
I_{31} \le C_* \|\theta\|_d^{\epsilon-a-p-s}  
 \int_0^{r_1} r^{p+a-1}
|\varrho(r)| dr,
\]
where $C_*=Ld_1^{p+s}\nu^{-p-s}
c_{p,d}$.
On the other hand for  $\|\theta\|_d \ge  r_1$, 
$\rho(\theta)=\varrho(\|\theta\|_d)$ is bounded from below as
\begin{align*}
\varrho(r_1)r_1^{-t_1}\|\theta\|_d^{t_1} \leq  \varrho(\|\theta\|_d).
\end{align*}
Therefore
\begin{align*}
I_{31} \le  \|\theta\|_d^{\epsilon-a-p-s - t_1} \times
C_*
\frac{r_1^{t_1}}{\varrho(r_1)} 
 \int_0^{r_1} r^{p+a-1}
|\varrho(r)| dr
\times \rho(\theta).
\end{align*}
Hence if $s> -t_1 - a - p$, then we can choose $\epsilon >0$ 
(say $ \epsilon=(a+p+s+t_1)/4$) such
that $ \epsilon-a-p-s-t_1<0$ and hence $ I_{31} \le C_{31}\rho(\theta)$ 
where
\[
 C_{31} =
C_*
 \frac{r_1^{t_1}}{\varrho(r_1)}\int_0^{r_1} r^{p+a-1}
|\varrho(r)| dr.
\]

Finally we consider $I_{32}$.  
Note $\varrho(r)\le \varrho(\|\theta\|_d)\|\theta\|_d^{-t_1} r^{t_1}$ 
for $r_1 \le r \le \|\theta\|_d$ and
\eqref{f-inequality}.
Then
\begin{align*}
I_{32} 
&\leq 
\|\theta\|_d^{\epsilon-a-t_1}
Ld_1^{p+s}\nu^{-p-s}\| \theta\|_d^{-p-s}
\varrho(\|\theta\|_d)
\int_{r_1 \leq \|x\|_d \leq \|\theta\|_d}  \|x\|_d^{t_1+a} dx
\\
&\leq \|\theta\|_d^{\epsilon-a-p-s-t_1} \times
C_*
\int_{r_1}^{\|\theta\|_d}  r^{p+a+t_1-1}dr \times \rho(\theta).
\end{align*}
Consider the integral $Q=\int_{r_1}^{\|\theta\|_d} r^{p+a+t_1-1}dr$.  If
$p+a+t_1 < 0$, then 
\[
Q \le r_1^{t_1 + a +p}/(-t_1-a-p).
\]
Therefore as in the case of $I_{31}$, 
if $s> -t_1 - a - p$, then we can choose $\epsilon>0$ 
(say $ \epsilon=(a+p+s+t_1)/4$) such
that $ \epsilon-a-p-s-t_1<0$ and hence 
\begin{align*}
 I_{32} \leq r_1^{\epsilon-s} C_{32}
 \rho(\theta) \leq C_{32}
 \rho(\theta)
\end{align*}
where $C_{32}=\{-1/(p+a+t_1)\}C_* $. 
If $ p+a+t_1 \geq 0$, 
\begin{align*}
 Q =\int_{r_1}^{\|\theta\|_d} r^{p+a+t_1+\epsilon_3-1-\epsilon_3}dr
\le
 \frac{\|\theta\|_d^{p+a+t_1+\epsilon_3}r_1^{-\epsilon_3}}
{p+a+t_1+\epsilon_3} 
\end{align*}
for any $ \epsilon_3>0$. 
Hence
\begin{align*}
 I_{32} \leq \frac{C_*}{p+a+t_1+\epsilon_3}
\|\theta\|_d^{\epsilon+\epsilon_3-s}\rho(\theta).
\end{align*}
If $ s>1$, we can choose $\epsilon$ and $ \epsilon_3$ 
(say $ \epsilon=\epsilon_3=1/4$)
such that
$ \epsilon+\epsilon_3-s<0$ and hence
\begin{align*}
 I_{32} \leq C_{32}\rho(\theta),
\end{align*}
where $C_{32}=(p+a+t_1+\epsilon_3)^{-1}C_*$.
%

We have now confirmed that if $s> \max(1, -t_1-a-p,a+t_2)$, 
there exist
$\epsilon>0$ and $C=C_1+C_2+C_{31}+C_{32}+C_{33}$, 
such that (\ref{eq:thm22}) folds for 
$\|\theta\|_d \ge 
\max(d_1 \nu^{-1} r_0,r_1, d_1(1-\nu)^{-1}r_1)$ which equals to
$ 2d_1\max(r_0,r_1)$ for $ \nu=1/2$.
\end{proof}

In the next section, we need asymptotic behavior of the expectation 
of $ \rho(X) \times h_i^\gamma(X)$ where $h_i(\theta)=H_i(\|\theta\|_d)$
given by \eqref{H} and $\gamma>0$.
\begin{corollary} \label{th2.2coro}
Assume {\bf F1}, {\bf B1} and {\bf B2}.  
For $a=0$ or $1$, $i\geq 1$, $\gamma>0$ and $j=1,\dots,p$, 
if $ 
s> \max(1, \gamma-t_1-a-p,t_2+a)$ and
$ \int_0^1 r^{p+a-1}| \varrho(r)|dr <\infty$,
there exists $ \epsilon>0$ (say $ \epsilon=\min(1,s+t_1+a+p-\gamma)
/4$) 
 such that
\begin{equation}
\label{eq:coro22}
\|\theta\|_d^{\epsilon-a}
\left|E_\theta[X_j^a \rho(X)h_i^\gamma(X)] - \theta_j^a 
\rho(\theta)h_i^\gamma(\theta)
 \right| 
< C \rho(\theta)h_i^\gamma(\theta)
\end{equation}
for  $\|\theta\|_d \geq 2d_1\max(r_0, r_1,\eta_0)$. 
Moreover $C$ does not depend on $i$.
\end{corollary}
\begin{proof}
Since 
we have
\[
 \eta\{\varrho(\eta)H_i^\gamma(\eta)\}'/\{\varrho(\eta)H_i^\gamma(\eta)\}
 = \eta 
 \varrho'(\eta)/\varrho(\eta)+\gamma\eta H'_i(\eta)/H_i(\eta),
\]
under Assumption {\bfseries B2} and by 5 of Theorem
 \ref{thm-beta-H},  for $ \epsilon_1= (s+t_1+a+p-\gamma)/16(>0) $
there exists $\eta_0$ such that
\begin{align*}
 t_1-\gamma-\gamma\epsilon_1 \leq 
    \eta \{\varrho(\eta)H_i^\gamma(\eta)\}'/\{\varrho(\eta)H_i^\gamma(\eta)\}
    \leq t_2
\end{align*}
for all $\eta \geq \max(\eta_0, r_1) $. 
Then for $ \epsilon= \min(1,s+t_1+a+p-\gamma)/4(>0) $,
 \eqref{eq:coro22} follows from Theorem \ref{thm-zenkin}.

For  $\lambda_1=\max(\eta_0, r_1) $,
\begin{align*}
 \frac{1}
{H_i^\gamma(\lambda_1)\varrho(\lambda_1)}\int_0^{\lambda_1}r^{p+a-1} 
|H_i^\gamma(r)\varrho(r)|dr \le
\frac{1}{H_1^\gamma(\lambda_1)\varrho(\lambda_1)}\int_0^{\lambda_1}r^{p+a-1} 
|\varrho(r)|dr
\end{align*}
which implies that $C$ does not depend on $i$.
\end{proof}

\section{Admissibility}
\label{sec:admissibility}
In this section, we give a sufficient condition for admissibility of the
generalized Bayes estimator with respect to a elliptically symmetric
prior density 
$g(\theta)=G(\| \theta \|_d)$.
The assumptions on $G$ are the following.

\begin{description}
\item[G1]  
$G(\eta)$ is continuously differentiable in $\eta>0$.
There exist $r_1 > 0$, $ t_1 \le t_2$ such
  that $ t_1\leq \eta G'(\eta)/G(\eta) \leq t_2$
 for all  $\eta \ge r_1$. 
\item[G1']  
By redefining $r_1$ if necessary in {\bfseries G1}, 
$t_1$ and $t_2$ can be taken as same sign.
\item[G1'']
$G'(\eta)$ is continuously differentiable in $\eta>0$.
there exist 
$ t_3 \le t_4$ such  that
$ t_3 \leq \eta G''(\eta)/G'(\eta) \leq t_4$
 for all  $\eta \ge r_1$. 
\item[G2] 
$\int_1^\infty \eta^{p-1}G(\eta)d\eta=\infty$ and there exists
$ 0<\gamma \leq 2$ such that 
$\int_1^\infty \eta^{p-1}G(\eta)H_1^\gamma(\eta)d\eta<\infty$.
\item[G3] $\lim_{\eta \to 0}\eta G'(\eta)/G(\eta)=t_0(>1-p). $
\item[FG1] $ \int_0^\infty r^{p-1}f(r)G(r)dr < \infty$ and $ \int_0^\infty
	   r^{p-2}F(r)G(r)dr < \infty$.
\end{description}

We discuss some implications of these assumptions.
By the assumption {\bfseries G3}, 
\begin{align*}
 \int_0^1 \eta^{p-1}G(\eta)d\eta<\infty \mbox{ and }
 \int_0^1 \eta^{p-1}|G'(\eta)| d\eta<\infty.
\end{align*}
 From the former
integrability and
{\bfseries G2}, the improperness of
   $g$ occurs only at infinity.
By 2 of Theorem \ref{thm-beta-H} and the assumption {\bf G2}, 
$ g(\theta)H_i^\gamma(\|\theta\|_d)$ for any fixed $i$ is integrable and hence
becomes a proper probability density by standardization.
Since $ H_i(\cdot)$ approaches $1$ as $i \to \infty$,
$ g(\theta)H_i^\gamma(\|\theta\|_d)$ 
is a sequence of proper densities approaching $g(\theta)$,
which is essential for using Blyth' method.

By {\bfseries G1} and 
Lemma \ref{lem:b3}, $G(\eta) = O(\eta^{t_2})$. Therefore if $t_2 < -p$,
then $g(\theta)=G(\|\theta\|_d)$ is a proper prior.  Since we are
considering an improper $g(\theta)$, we assume $t_2 \ge -p$ from now on.
Moreover $\eta^{t_1}/G(\eta)=O(1)$ by the assumption {\bfseries G1}.
Since $H_1 \in \mbox{RV}_{-1}$, $ \eta^{p-1}\eta^{t_1}H_1^2(\eta) $ for
$ t_1>2-p$ is not integrable
at infinity
and $ \eta^{p-1}G(\eta)H_1^2(\eta)$ is  not so either.
Hence we also assume $t_1 \leq 2-p$.
We now discuss when 
the integrability of $\int_1^\infty
\eta^{p-1}G(\eta)H_1^\gamma(\eta)d\eta $ in {\bfseries G2} holds
and the relationship with {\bfseries G1}.
When we take $ \beta(\eta)$ as in \eqref{beta1}, 
we easily see that there exists $L_1$ such that
$ H_1(\eta)\leq L_1 (\eta+c)^{-1} \{\log(\eta+c)\}^{-1}$.
If $t_2 \leq 2-p$, there clearly exists $L_2$ such that
$ G(\eta) \leq L_2\eta^{2-p}$ for $ \eta>1$.
Hence
\begin{align*}
\int_1^\infty \eta^{p-1}G(\eta)H^2_1(\eta)
d\eta & \leq 
 L_1^2 L_2 
 \int_1^\infty (\eta+c)^{-1}  \{\log(\eta+c)\}^{-2} d\eta \\
& = L_1^2 L_2  \{\log(1+c)\}^{-1}<
\infty ,
\end{align*}
which shows that there exists $\gamma(=2)$ for the integrability in 
{\bfseries G2}.
On the other hand, if $ t_2>2-p$ then the integrability in 
{\bfseries G2} may not be apparent. But the condition 
$t_2 \le 2-p$ is not a necessary
condition for the integrability.
%

%
If $G(\eta) $ is regularly varying, then for any $\epsilon>0$, we can
choose $r_1$, $t_1, t_2$, $t_3, t_4$,  such that $t_2-t_1<\epsilon $ and $
t_4-t_3<\epsilon$ for all $\eta \geq r_1$. 
However in ${\bf G1}$ we are allowing the case that
$\liminf_{\eta\rightarrow\infty} \eta G'(\eta)/G(\eta)$ is strictly less
than $\limsup_{\eta\rightarrow\infty} \eta G'(\eta)/G(\eta)$.  Hence we
are dealing with a broader class of $G(\eta)$ than the class of regularly
varying functions. It should also be noted that $ t_4 $ and $t_3$ are not
always smaller than $ t_2$ and $t_1$, respectively. See
\cite{geluk-dehaan1987} for the detail.

The generalized Bayes estimator $ \delta_{g}$ 
with respect to the improper density 
$g(\theta)$ is written as
\begin{align}
 \delta_g (x)&=\frac{\int_{R^p} \theta f(\|x - \theta\|)g(\theta)d\theta}
{\int_{R^p}  f(\|x - \theta\|)g(\theta)d\theta} 
\nonumber \\
&=x+ \frac{\int_{R^p} (\theta-x) f(\|x - \theta\|)g(\theta)d\theta}
{\int_{R^p}  f(\|x - \theta\|)g(\theta)d\theta} 
\nonumber \\
&=x+ \frac{\int_{R^p}  F(\|x - \theta\|) \nabla
 g(\theta)d\theta}
{\int_{R^p}  f(\|x - \theta\|)g(\theta)d\theta} \label{bayes}, 
\end{align}
which is well-defined if both $ \int_{R^p} F(\|x - \theta\|) \nabla
g(\theta)d\theta $ and $ \int_{R^p} f(\|x - \theta\|)g(\theta)d\theta $ are
integrable for all $x$.  These are guaranteed by the assumption {\bfseries
FG1} and Lemma \ref{thm-zenkin-1} in the below.  Write
\begin{align*}
 m(\psi | x) &= \int_{R^p} \psi(\theta) f(\|\theta - x\|)d\theta \\
M(\psi | x) &=  \frac{1}{C_f}\int_{R^p} \psi(\theta) F(\| \theta - x \|) 
 d\theta
\end{align*}
where $ C_f=\{\pi^{p/2}/\Gamma(p/2+1)\} \int_0^\infty z^{p+1}f(z)dz $.
%
%
Notice that $F( \cdot )/C_f$ is 
 a probability density function because
\begin{align*}
\int_{R^p} \|y - \theta \|^\alpha F(\|y - \theta \|)dy
& = \int_{R^p} \|y\|^\alpha \{ \int_{\| y  \|}s f(s)ds \} dy \\
&= c_p \int_0^{\infty} r^{p-1+\alpha} \int_{r}^\infty s f(s)ds dr \\
&= c_p \int_0^\infty r^{p+1+\alpha} \int_1^\infty t f(r t)dt dr \\
& = c_p \int_1^\infty t \{ \int_0^{\infty } r^{p+1+\alpha} f(r t ) dr \} dt \\
& = c_p \int_1^\infty t^{-p-1-\alpha}dt \cdot \int_0^\infty z^{p+1+\alpha}f(z)dz \\
& =\frac{c_p}{p+\alpha} \int_0^\infty z^{p+1+\alpha}f(z)dz.
\end{align*}
Then $\delta_g$ is written as
\[
\delta_g (x)= x + C_f \frac{M(\nabla g|x)}{m(g|x)}.
\]
Note that by {\bf G1} the $j$-th element of $\nabla g$ is given by
\[
\nabla_j g(\theta)= d_j^2 \theta_j \frac{G'(\|\theta\|_d)}{\|\theta\|_d}.
\]
We also write
\[
h_i(x)=H_i(\|x \|_d).
\]

Now we state the following lemma in preparation of our main theorem.

\begin{lemma} 
\label{thm-zenkin-1}
\begin{enumerate}
\item 
Assume {\bf G1}, {\bfseries G3} and {\bf F1}.
Then 
\begin{align}
 m( g|x) & \approx g(x) & \mbox{ if }s &> \max(1, -t_1-p, t_2) \label{y1}\\
 m( gh_i^\gamma|x) & \approx g(x)h_i^\gamma(x) & \mbox{ if }s&> \max(1,
 \gamma-t_1-p, t_2) \label{y2}\\ 
 M( gh_i^\gamma|x) & \approx g(x)h_i^\gamma(x) & \mbox{ if }s &> 2+\max(1,
 \gamma-t_1-p, t_2) \label{y3}\\ 
 M( g\|\theta\|^{-1}_d|x) & \approx g(x)\|x\|^{-1}_d & \mbox{ if }s &> 2+\max(1,
 1-t_1-p, t_2-1) \label{y4} \\
 M( gh_i^\gamma \|\theta\|^{-1}_d|x) & \approx g(x)h_i^\gamma(x)
\|x\|^{-1}_d & \mbox{ if }s &> 2+\max(1,
 \gamma+1-t_1-p, t_2-1) \label{y5}. 
\end{align}
\item Assume {\bf G1'}, {\bfseries G1''}, {\bfseries G3} and {\bf F1}.
Then
\begin{align}
M(\nabla_j g |x) & \approx
\nabla_j g (x)  & \mbox{ if }s &> 2+\max(1, -t_3-p, t_4) \label{y6}\\
M(\nabla_j gh_i^\gamma |x) & \approx
\nabla_j g (x)h_i^\gamma(x)  & \mbox{ if }s &> 2+\max(1, \gamma-t_3-p,
 t_4).  \label{y7}
\end{align}
\item $ M(\nabla g|x)/m(g|x)$ is bounded in $x$ if
$s>2+\max(1,1-t_1-p,t_2-1)$.
\end{enumerate}
\end{lemma}

\begin{proof}
When we consider the asymptotic behavior of $ M(\cdot|x)$, that is,
the expectation under the probability density $F(\|\theta - x\|)/C_f$,
we have only to substitute $s$ for $s-2$ in order to have corresponding 
results for Theorem \ref{thm-zenkin} and Corollary \ref{th2.2coro}
because under the assumption {\bfseries F1} there exists $L_1$ such
that $r^{p+s-2}F(r)<L_1$ for all $r\ge r_0$.
We easily see that
\eqref{y1}, \eqref{y4} and \eqref{y6} follow from 
Theorem \ref{thm-zenkin} and that \eqref{y2}, \eqref{y3}, \eqref{y5} 
and \eqref{y7} follow from
 Corollary \ref{th2.2coro}.

Note that by the assumptions {\bfseries G1} and {\bfseries G3}
there exists $L_2$ such that $ |\eta G'(\eta)/G(\eta)|<L_2$
 for all $ \eta>0$.
Then $ M(\nabla g|x)/m(g|x) \leq d_1 M(|G'| \, |x)/m(g|x)
\leq L_2d_1 M(G/ \| \theta\| \, |x)/m(g|x) $.
The value at $x=0$ is clearly bounded under the assumption {\bfseries
 FG1} and the value at $\|x\|_d \to
\infty$  is bounded by \eqref{y1} and
\eqref{y4} if $s> 2+\max(1, 1-t_1-p, t_2-1)$.
\end{proof}

By part 3 of Lemma \ref{thm-zenkin-1}, 
$ M(\nabla g|x)/m(g|x)$ is bounded in $x$ and hence
the risk function of $\delta_g$
is finite  because 
\begin{align*}
 R(\theta,\delta_g)&=E[\|X-\theta+C_f M(\nabla g|X)/m(g|X)\|_Q^2]\\
&\leq Q_{\mathrm{max}} E[\|X-\theta+C_f M(\nabla g|X)/m(g|X)\|^2]\\
&\leq 2 Q_{\mathrm{max}} \, \{ E[\|X-\theta\|^2]+C_f^2 E[\|M(\nabla
 g|X)/m(g|X)\|^2]\},
\end{align*}
where $Q_{\mathrm{max}}$ is the largest eigenvalue of $ Q$.

Now we state the main theorem of this paper.
\begin{thm}
\label{thm:main}
\begin{enumerate}
 \item 
Assume {\bf G1}, {\bfseries G2}, {\bfseries G3}, {\bf F1}, {\bf FG1}.
Then 
the generalized Bayes estimator with respect to $g$ is admissible
 if 
$ s> 2+\max(1, 1+\gamma-t_1-p)$ and $t_2<2-p$.
\item
We also assume {\bfseries G1'} and {\bfseries G1''}. 
Then 
the generalized Bayes estimator with respect to $g$ is admissible
 if 
$ s> 2+\max(1, 3-t_1-p, t_2, 2-t_3-p, t_4)$ and 
$t_2 \ge 2-p$.
\end{enumerate}
\end{thm}
Although the moment conditions for $s$ in the theorem above 
looks complicated, 
it is just from the assumptions {\bfseries G1}, {\bfseries G1'} and
{\bfseries G1''} which make our class of
$G$ broader than the class of regularly varying functions.
We see that
the condition reduces to $s>3$ for regularly varying functions $G$.
Before giving a proof of the main theorem we present it as a 
corollary.
\begin{corollary} \label{main-coro-2}
Suppose that $G(\eta)$ is regularly varying with index $k$ for
$-p \le k \le 2-p$.
Assume {\bfseries G3}, {\bfseries F1} with $ s> 3$ and {\bfseries FG1}.
\begin{enumerate}
 \item Assume  $-p \le k < 2-p$ and 
$ G(\eta)$ is continuously differentiable.
  Then the generalized Bayes estimator with respect to
$g$ is admissible.
\item Assume $k=2-p$ and $ G(\eta)$ satisfies 
\begin{align}
\textstyle G(\eta) \leq  \eta^{2-p} \{\int_\eta^\infty
\beta(r)dr\}^2\{\eta\beta(\eta)\}^{-1} \ \mbox{ for } \eta \geq 1, \label{girigiri}
\end{align}
 and 
 is twice continuously differentiable. 
  Then the generalized Bayes estimator with respect to
$g$ is admissible.
\end{enumerate}
\end{corollary}
\begin{proof}
When $ -p \le k < 2-p $, 
we can take $t_1=k-\epsilon/3$ and $ t_2=k+\epsilon/3$ for 
any $ 0<\epsilon<\{3/2\}(2-p-k)$ in {\bfseries G1}. 
Let $\gamma=k+p+2\epsilon/3$. Clearly 
$ 0<\gamma<2$ and it satisfies the integrability in {\bfseries G2}.
Since $ t_2<2-p$, we have only to apply part 1 of Theorem \ref{thm:main}.
$  1+\gamma-t_1-p = 1+\epsilon$ and hence
the moment condition for $ k<2-p$ is $s>3$.

Next we consider the case $k=2-p$. 
Note that there exists $ L_1$ such that $ H_1(\eta)\leq L_1 \beta(\eta)/
{\int_{\eta}^\infty \beta(r)dr}$ by part 2 of Lemma
2.1. So if $G(\eta)$ satisfies
\eqref{girigiri}, we have
\begin{align*}
\int_1^\infty \eta^{p-1}G(\eta)H_1^2(\eta)d\eta
 & \leq 
L_1 \int_1^\infty \beta(\eta)d \eta<\infty,
%
\end{align*}
which shows that the integrability in the assumption {\bfseries G2}
is guaranteed.
We can take $t_1=2-p-\epsilon$, $ t_2=2-p+\epsilon$,
$t_3=1-p-\epsilon$ and $t_4=1-p+\epsilon$ for any $ \epsilon>0$
in the assumptions
 {\bfseries G1}, {\bfseries G1''}.
Since $ t_2 \geq 2-p$, we have only to apply part 2 
of Theorem \ref{thm:main}.
Since
\begin{align*}
 3-t_1-p=1+\epsilon, \quad t_2 \leq 1+\epsilon, \quad 2-t_3-p=1+\epsilon,
\quad t_4 \leq \epsilon
\end{align*}
the moment condition for $ k=2-p$ is $s>3$.
\end{proof}
In particular the boundary case in \eqref{girigiri} 
by taking $\beta(\eta)$ in
(\ref{beta1}) was the motivating one for this paper.
\begin{corollary} \label{main-coro}
Assume {\bfseries F1} with $ s> 3$ and {\bfseries FG1}.
  Then the generalized Bayes estimator with respect to
$ \|\theta\|_d^{2-p} \prod_{i=0}^n{\rm Log}_i(\|\theta\|_d+c)$,
where $n$ is a nonnegative integer and ${\rm Log}_n(c)>0$, 
is admissible.
\end{corollary}
In the normal case, \cite{Brown1971}'s sufficient conditions for
admissibility and inadmissibility are known.  He showed that the
generalized Bayes estimator with respect to $g$ is admissible if
\begin{align}
\label{brown-condition}
 \int_{\|x\|_d>1} \|x\|_d^{2-2p}\{m(g|x)\}^{-1} dx
\end{align}
diverges and inadmissible if \eqref{brown-condition} converges.
By Lemma \ref{thm-zenkin-1}, we see that 
\begin{align*}
 (1/2)G(\|x\|_d)<m(g|x)<2G(\|x\|_d)
\end{align*}
for sufficiently large $\|x\|_d$
and hence that $G(\eta) \leq \eta^{2-p}\prod_{i=0}^n{\rm Log}_i(\eta+c)
$ leads to admissibility and 
$G(\eta) \geq \eta^{2-p}\prod_{i=0}^{n-1}{\rm Log}_i(\eta+c)
 {\rm Log}_n^2(\eta+c)$ leads to inadmissibility.
Therefore our sufficient condition in Theorem
\ref{thm:main} is very close to being necessary.

We also notice that the prior density suggested in Corollary 
\ref{main-coro} becomes $ |\theta| \prod_{i=0}^n{\rm Log}_i(|\theta|+c)$
and $ \prod_{i=0}^n{\rm Log}_i(\|\theta\|_d+c)$ for $p=1, \ 2$
respectively, which are thicker than the Lebesgue measure.
In the normal case, \cite{Brown1971} has already pointed it out.

Furthermore we indicate that our moment condition $s>3$ is very tight
because, as pointed out in \cite{Perng1970}, admissibility requires the
existence of moment one degree higher than what is needed for finite risk in 
various estimation problems.

Now we give a proof of Theorem \ref{thm:main}.

\begin{proof}[Proof of Theorem \ref{thm:main}]
Let $\delta_{gi}$ denote the Bayes estimator with respect to the proper
prior density $g(\theta)h_i^\gamma(\theta)$. Then the Bayes risk difference
of $ \delta_g$ and $ \delta_{gi}$ with respect to the density 
$ g(\theta)h_i^\gamma(\theta)$ is written as
\begin{align*}
 \Delta &= \int_{R^{p}} \left[
R(\theta,\delta_{g})-R(\theta,\delta_{gi}) \right]
 g( \theta)h_i^\gamma(\theta )d \theta \\
&= \int_{R^p} \int_{R^p} [ \| \delta_{g}- \theta \|_Q^{2}- \| \delta_{gi}- \theta
 \|_Q^{2}]f(\|x-\theta\|)g(\theta)h_i^\gamma(\theta)
 d\theta dx\\
&= \int_{R^p} \bigg\{ [ \| \delta_{g} \|_Q^{2}-\| \delta_{gi}\|_Q^{2}]
  \int_{R^p} f(\| x- \theta\|) g(\theta)h_i^\gamma(\theta)d\theta
\\ & \quad 
-2(\delta_{g}-\delta_{gi})Q'\int_{R^p} \theta f(\|x-\theta\|)
g(\theta)h_i^\gamma(\theta)
 d\theta \bigg\}dx \\
&= \int_{R^{p}} \| \delta_{g} - \delta_{gi} \|_Q^{2} \left\{\int_{R^p} 
 f(\|x-\theta\|)g(\theta)h_i^\gamma(\theta)
 d\theta\right\} dx \\
&= C_f^2\int_{R^{p}} \left\| \frac{M(\nabla g|x)}{m(g|x)} - 
\frac{M(\nabla \{gh_i^\gamma\}|x)}{m(gh_i^\gamma|x)}
\right\|_Q^2 m(gh_i^\gamma|x)dx \\
&= C_f^2\int_{R^{p}} \left\| \frac{M(\nabla g|x)}{m(g|x)} - 
\frac{M(\nabla g h_i^\gamma|x)}{m(gh_i^\gamma|x)}-\frac{M(g \nabla h_i^\gamma|x)}{m(gh_i^\gamma|x)}
\right\|_Q^2 m(gh_i^\gamma|x)dx .
\end{align*}

In the same way as in \cite{Brown-Hwang1983},
we have
\begin{align*}
 \Delta & \leq 2C_f^2 Q_{\mathrm{max}}
 \int_{R^{p}} \left\| \frac{M(\nabla g|x)}{m(g|x)} - 
\frac{M(\nabla g h_i^\gamma|x)}{m(gh_i^\gamma|x)}
\right\|^2 m(gh_i^\gamma|x)dx \\
& \qquad + 2 C_f^2Q_{\mathrm{max}}
\int_{R^{p}} \left\| \frac{M(g \nabla h_i^\gamma|x)}{m(gh_i^\gamma|x)}
\right\|^2 m(gh_i^\gamma|x)dx  \\
&=2C_f^2Q_{\mathrm{max}}(B_i+A_i). \qquad  \mathrm{(say)}
\end{align*}

Using the Cauchy-Schwartz inequality for $A_i$, we have
\begin{align*}
 A_i &=\gamma^2 \int_{R^{p}} \left\| M(g h_i^{\gamma-1} \nabla h_i|x)
\right\|^2 \{ m(gh_i^\gamma|x) \}^{-1} dx \\
& \leq \gamma^2 \int_{R^{p}} \frac{M(gh_i^\gamma|x)}{m(gh_i^\gamma|x)}
 M(g h_i^{\gamma-2} \| \nabla h_i \|^2|x)
 dx \\
& \leq \gamma^2 \int_{R^{p}} \frac{M(gh_i^\gamma|x)}{m(gh_i^\gamma|x)}
 M(g h_1^{\gamma-2} \| \nabla h_i \|^2|x)
 dx
\end{align*}
for $ 0<\gamma\leq 2$.
The ratio $M(gh_i^\gamma|x)/m(gh_i^\gamma|x)$ is bounded from above
by $M(g|x)/m(gh_1^\gamma|x)$ and hence the value at $x=0$ is clearly
bounded under the assumption {\bfseries FG1}.
By \eqref{y2} and \eqref{y3}, we have
\begin{align*}
\lim_{\|x\|_d \to \infty} \frac{M(gh_i^\gamma|x)}{m(gh_i^\gamma|x)}=1
\end{align*}
uniformly in $i$.  This implies that there exists $c_1$ such that
$M(gh_i^\gamma|x)/m(gh_i^\gamma|x)<c_1$ for all $x$ and for all $i$. Then
\begin{align*}
A_i & \leq \gamma^2c_1\int_{R^{p}}  M(g h_1^{\gamma-2} 
\| \nabla h_i  \|^2 |x) dx \\
& = \gamma^2 c_1 \int_{R^{p}}\{1/C_f\} F(\| x - \theta \|)dx 
\int_{R^{p}} g(\theta) 
h_1^{\gamma-2}(\theta) \| \nabla h_i (\theta)\|^2
 d\theta \\
&= \gamma^2 c_1 \int_{R^{p}} g(\theta) h_1^{\gamma-2}(\theta) 
\| \nabla h_i(\theta) \|^2 d\theta.
\end{align*}
By 4 of Theorem \ref{thm-beta-H} we have $  \| \nabla h_i(\theta) \|
< 2 d_1 \beta(\|\theta\|_d)/\int_{\|\theta\|_d}^\infty \beta(r)dr $
and together with 2 of Theorem \ref{thm-beta-H}
$ \| \nabla h_i(\theta) \|/h_1(\theta)$ for all $\theta$
is bounded from above by $L_3$ independent of $i$. Therefore
\[
A_i \leq 
\gamma^2 c_1L_3^2 \int_{R^{p}} g(\theta) h_1^{\gamma}(\theta) 
d\theta,
\]
which is bounded by the assumption {\bfseries G2}.
Furthermore $ \| \nabla h_i(\theta)\| \to 0$ as $i \to \infty $ by 3 of
Theorem \ref{thm-beta-H}. Therefore by the dominated convergence theorem $
A_i $ converges to $0 $ as $i \to \infty$.

Next we consider $ B_i$.
$ M(\nabla g|x)$ and $M(\nabla g h_i^\gamma|x)$ at $x=0$ are zero vectors 
because $g$ and $h_i^\gamma$ are function of $ \| \theta \|_d$.
So the integrand of $B_i$ is bounded around $x=0$.
For the asymptotic property of the integrand of $B_i$,
we need to distinguish two cases: $t_2<2-p$ and $t_2\geq 2-p$.
When $t_2<2-p$, we can bound the norm in the integrand of $B_i$
from above somewhat roughly.
Using \eqref{y1}, \eqref{y2}, \eqref{y4} and \eqref{y5} in Lemma
 \ref{thm-zenkin-1} and noting
that $ |\eta G'(\eta)/G(\eta)|<L_2$ by the assumptions 
{\bfseries G1} and {\bfseries G3}
we have
\begin{align*}
 \frac{1}{d_j^2}\left| \frac{ M(\nabla_j g|x)}{m(g|x)}-\frac{ M(\nabla_j g
 h_i^\gamma|x)}{m(gh_i^\gamma|x)} \right| 
&  <  \frac{ M(|G'|\, |x)}{m(g|x)}
+\frac{ M(|G'|h_i^\gamma \, |x)}{m(gh_i^\gamma|x)}
\\
&  <  L_2 \left(\frac{ M(G \|\theta\|^{-1}\, |x)}{m(g|x)}
+\frac{ M(G h_i^\gamma\|\theta\|^{-1} |x)}{m(gh_i^\gamma|x)} \right)\\
&  < c\|x\|_d^{-1}
\end{align*}
for all sufficiently large $\|x\|_d$ and for all $i$.
When $t_2 \geq 2-p$, we have to bound it from above more strictly.
By \eqref{y1}, \eqref{y2}, \eqref{y6} and \eqref{y7} in
Lemma \ref{thm-zenkin-1},  we have
\begin{align*}
& \frac{1}{d_j^2}\left| \frac{ M(\nabla_j g|x)}{m(g|x)}-\frac{ M(\nabla_j g
 h_i^\gamma|x)}{m(gh_i^\gamma|x)} \right| \\
&\qquad  = \left| \frac{G'(\|x\|_d)}{G(\|x\|_d)\|x\|_d}
\frac{x_j+O(\|x\|_d^{1-\epsilon_0}) }{1+O(\|x\|_d^{-\epsilon_0})} -
\frac{G'(\|x\|_d)}{G(\|x\|_d)\|x\|_d}
\frac{x_j+O(\|x\|_d^{1-\epsilon_0}) }{1+O(\|x\|_d^{-\epsilon_0})}
\right| \\
&\qquad  < c\|x\|_d^{-1-\epsilon_0} \\
&\qquad  < H_1(\|x\|_d) 
\end{align*}
for some $\epsilon_0$,
for all sufficiently large $\|x\|_d$ and for all $i$.
Moreover $m(gh_i^\gamma|x)\leq m(g|x)$ and $m(g|x) < 2G(\|x\|_d)$
for all sufficiently large $\|x\|_d$ by \eqref{y1}.  
Therefore there exist
$C_1$, $C_2$, $C_3$ and $C_4$ such that the integrand of $B_i$ is less
than
\begin{align*}
\begin{cases}
\min\{C_1, C_2 \|x\|_d^{-2+t_2}\} & \mbox{ if } t_2<2-p \\
\min\{C_3, C_4 G(\|x\|_d)H_1^2(\|x\|_d)\} & \mbox{ if } t_2 \ge 2-p.
\end{cases}
\end{align*}
Therefore 
$ B_i $ converges to $0 $ as $i \to \infty$ by the dominated convergence
theorem.

Finally we confirm that we use \eqref{y1}--\eqref{y5} for $t_2<2-p$
and \eqref{y1}--\eqref{y3}, \eqref{y6} and \eqref{y7} for $ t_2 \geq
 2-p$. Note also $\max(1, t_2)=1$ for $t_2<2-p$ in the moment condition.
\end{proof}

\section{The generalized Bayes estimator with respect to
 the harmonic  prior and its minimaxity} 
In this section, we show that the generalized Bayes estimator with
respect to the harmonic prior has a form simple enough to check 
some sufficient conditions for minimaxity under the quadratic
loss function $L_I(\theta,\delta)=\| \delta- \theta \|^2$
given in early studies.
We demonstrate that it is minimax for some $f$.

In \eqref{bayes}, the generalized Bayes estimator can be also written as
\begin{align*}
 \delta_g(x)=x+C_f \frac{\nabla_x M(g|x)}{m(g|x)}.
\end{align*}
For $p \geq 3$ and $g(\theta)=\|\theta\|^{2-p}$, 
we have
\begin{align}
  m(g|x)&=\int_{R^{p}}f(\|x-\theta\|)\| \theta \|^{2-p}d\theta 
 =\int_{R^{p}}f(\|\eta\|)\| x-\eta \|^{2-p}d\eta \nonumber \\
 & \ \ = 
c_{p-1}\int_{0}^{\infty}\int_{0}^{\pi}f(\lambda)(\lambda^2
+2\lambda r\cos \varphi +r^2)^{1-p/2}\lambda^{p-1}
\sin^{p-2}\varphi d\lambda d \varphi \nonumber \\
 & \ \ = 
c_{p-1} r^2 \int_{0}^{\infty}\int_{0}^{\pi}f(rt)
(1+2t\cos \varphi+t^2)^{1-p/2}t^{p-1}\sin^{p-2}\varphi 
dt d \varphi \nonumber \\
 & \ \ = 
c_{p}\left(r^2\int_{0}^{1}t^{p-1}f(rt)dt+r^2\int_{1}^{\infty}
t f(rt)dt \right)
 \nonumber \\ 
& \ \ = 
c_{p}\left( \left[ -t^{p-2} F(rt) \right]_0^1+
(p-2)\int_0^1 t^{p-3}F(rt)dt+ F(r)
 \right)
 \nonumber \\ 
 & \ \ = 
c_p(p-2)\int_{0}^{1}t^{p-3}F(rt)dt,
\end{align}
where $ r=\|x\| $.
The fifth equality in the above equation follows from
the relation
\begin{equation*}
 \int_{0}^{\pi}(1+2t\cos \varphi +t^2)^{1-p/2}
\sin^{p-2}\varphi d\varphi=B(p/2-1/2,1/2)\min(t^{2-p},1),
\end{equation*}
which is proved in Lemma \ref{int} in the end of this section.
In the same way, we have
\begin{align*}
C_f \nabla_x M(g|x)=  -x c_p(p-2)\int_{0}^{1}t^{p-1}F(rt)dt.
\end{align*}
Hence 
the generalized Bayes estimator is written as
$ \delta_{*}(X)=(1-\phi_{*}(\|X\|)/\|X\|^2)X $,  where 
\begin{align*}
 \phi_{*}(r)=r^2\frac{\int_{0}^{1}t^{p-1}F(rt)dt}
{\int_{0}^{1}t^{p-3}F(rt)dt}.
\end{align*}
Some properties of the behavior of $ \phi_{*}(r) $ are easily derived as
follows. 
\begin{thm} \label{lemma1}
 \begin{enumerate}
  \item $ \lim_{r \rightarrow \infty }\phi_{*}(r)
=(p-2)E_{0}(\|X\|^{2})/p $.
\item $ \phi_{*}(r) $ is  nondecreasing in $r $ for any $ f$. 
\item $ \phi_{*}(r)/r^2 $ is nonincreasing in $r $ 
if $ F(t)\{t^2f(t)\}^{-1} $ is nonincreasing.   
 \end{enumerate}
\end{thm}
\begin{proof}
$ \phi_{*}(r) $ can be written as
$  \int_{0}^{r}t^{p-1}F(t)dt/
\int_{0}^{r}t^{p-3}F(t)dt$
and we have
\begin{equation*}
 \lim_{r \rightarrow \infty }\phi_{*}(r)= 
\frac{\int_{0}^{\infty}t^{p-1}F(t)dt}
{\int_{0}^{\infty}t^{p-3}F(t)dt} 
=\frac{p-2}{p}\frac{\int_{0}^{\infty}t^{p+1}f(t)dt}
{\int_{0}^{\infty}t^{p-1}f(t)dt}=\frac{p-2}{p}E_{0}[\|X\|^{2}]. 
\end{equation*}
The derivative of $ \phi_{*}(r) $ is calculated as
\begin{align*}
 \phi'_{*}(r)=\frac{r^{p-3}F(r)}
{\left( \int_{0}^{r}t^{p-3}F(t)dt \right)^2}
\int_{0}^{r}(r^2-t^2)t^{p-3}F(t)dt,
\end{align*}
which is nonnegative for any $f$. The derivative of $ \phi_{*}(r)/r^2 $ is
 calculated as
\begin{align*}
 \frac{d}{dr}\left(\phi_{*}(r)/r^2 \right) & =
r \left( \int_{0}^{1}t^{p-3}F(rt)dt \right)^{-2} \left( 
 \int_{0}^{1}t^{p-1}F(rt)dt\int_{0}^{1}t^{p-1}f(rt)dt \right. \\
 & \qquad - \left.
\int_{0}^{1}t^{p+1}f(rt)dt\int_{0}^{1}t^{p-3}F(rt)dt \right).
\end{align*}
If $ F(t)\{t^2f(t)\}^{-1} $ is nonincreasing, 
the right-hand side of the equality above is nonpositive by the
 covariance  inequality.  
\end{proof}

Now we consider the minimaxity of $\delta_*$. 
We present a brief list of known sufficient conditions for minimaxity given
in previous papers, for the estimator of the form \eqref{orthogonal} with
nonnegative and nondecreasing $\phi(r)$.

\medskip
\footnotesize
\begin{center}
\begin{tabular}{|c|c|c|c|}\hline
Author & $ p $ & $ \phi(r)/r^2 $ & upper bound of $\phi$\\ \hline \hline
\multicolumn{4}{|c|}{general} \\ \hline
\cite{Berger1975} & $ p \geq 3 $ & &  
$ 2(p-2)\inf_{s \in U}F(s)/f(s)  $   \\ \hline
\cite{Brandwein1979} & $ p \geq 4 $ & $ \searrow $ & 
$ 2(p-2)(pE_{0}(\|X\|^{-2}))^{-1} $ \\ \hline \hline
 \multicolumn{4}{|c|}{unimodal or $ f $ is nonincreasing} \\ \hline
\cite{Brand-Straw1978} & $ p \geq 4 $ & $ \searrow $ &  
$ 2p((p+2)E_{0}(\|X\|^{-2}))^{-1} $ \\ \hline
\cite{Ralescu-etal1992} & $p=3 $ & $ \searrow $ &
$ 0.93(E_{0}(\|X\|^{-2}))^{-1} $ \\ \hline \hline
 \multicolumn{4}{|c|}{$ F(t)/f(t) $ is nondecreasing} \\ \hline
\cite{Bock1985}  & $ p \geq 4 $ & $ \searrow $ & 
$ 2(E_{0}(\|X\|^{-2}))^{-1}$ \\ \hline
\hline
 \multicolumn{4}{|c|}{scale mixtures of multivariate normal} \\ \hline
\cite{Strawderman74}   & $ p \geq 3 $ & $ \searrow $ & 
$ 2(E_{0}(\|X\|^{-2}))^{-1}$ \\ \hline
\end{tabular}
\end{center}
\normalsize
\medskip
In the table, $ U=\{t \geq 0| f(t)>0 \}$ and an arrow $ \searrow $ means {\it
nonincreasing}. 
It is noted that 
$ f(t) $ is nonincreasing in $t $ if $ F(t)/f(t) $ is nondecreasing in 
$t $ and that $ F(t)/f(t) $ is nondecreasing in $t $ if $f$ is a scale 
mixtures of multivariate normal.

Combining Theorem \ref{lemma1} and  the table above, 
we can derive a sufficient condition for minimaxity of $ \delta_{*}(X) $
and we state it in the following theorem for $p \geq 4$.

\begin{thm}
\label{thm:minimaxity}
\begin{enumerate}
\item Assume $ t^{-2}F(t)/f(t)$ is nonincreasing.
\begin{enumerate}
\item $\delta_*$ is minimax if $ E_0[\|X\|^2]E_0[\|X\|^{-2}]
\leq 2$.
\item Assume also $ f(t)$ is nonincreasing.
Then $\delta_*$ is minimax \\ if $ (p^2-4)E_0[\|X\|^2]E_0[\|X\|^{-2}]
\leq 2p^2$.
 \item Assume also $ F(t)/f(t)$ is nondecreasing.
Then $\delta_*$ is minimax \\ if $ (p-2)E_0[\|X\|^2]E_0[\|X\|^{-2}]
\leq 2p$.
\end{enumerate}
\item Assume $0<\inf_{s \in U}F(s)/f(s)<\infty$.
Then $\delta_*$ is minimax \\ if $ E_0[\|X\|^2]\leq 2p\inf_{s \in
 U}F(s)/f(s)$. 
\end{enumerate}
\end{thm}
\cite{Berger1975} and \cite{Bock1985} gave several examples
of $f$, checked the monotonicity of $f(t)$, $F(t)/f(t)$, and
$t^{-2}F(t)/f(t)$ and  calculated an upper bound of $\phi(r)$.
In this paper we give just two examples but we believe that
the estimator $ \delta_{*}(X) $ is  minimax for a broad class of 
spherically symmetric distributions.
\begin{example}
We consider $  f(s)=s^{\alpha}\exp(-\beta s^2) $ for $ \alpha,\ \beta>0 $. 
We have
\begin{equation*}
 \frac{F(t)}{t^2f(t)}=\int_{1}^{\infty}u^{\alpha+1}\exp(\beta t^2 (1-u))du,
\end{equation*}
which is decreasing in $ t $. By an integration by parts, we have
\begin{align*}
 \frac{F(t)}{f(t)}&=\frac{1}{2\beta}+\frac{\alpha}{2\beta}
\frac{\int_t^\infty s^{\alpha-1}\exp(-\beta s^2)ds}{t^\alpha
\exp(-\beta t^2)}\\
&=\frac{1}{2\beta}+\frac{\alpha}{2\beta}\int_1^\infty
u^{\alpha-1}\exp(\beta t^2\{ 1-u^2\})du
\end{align*}
and hence 
$ \inf_{t \geq 0}F(t)/f(t)=(2\beta)^{-1} $.
We also have $ E_{0}(\|X\|^{2})
=(p/2+\alpha/2)/\beta $ and $E_{0}(\|X\|^{-2})^{-1}=(p/2+\alpha/2-1)/\beta $.
Therefore the generalized Bayes estimator is minimax 
if $ \alpha \leq p $ for $ p \geq 3 $ by \cite{Berger1975}'s conditions 
and if $ \alpha \geq 4-p $ for $ p \geq 4 $ by \cite{Brandwein1979}'s
 conditions regardless of $ \beta $. 
Hence the estimator for $ p \geq 4 $ is minimax regardless of $ \alpha $
 and $ \beta$.
\end{example}

%

\begin{example}
We consider 
$  f(t)=\exp(-t^2/2)-a\exp(-t^2/\{2b\}) $ for $ 0< a \leq 1 $, $ 0<b<1 $. 
Note that if $ a \leq b $ then $f $ is unimodal and if $ a > b $ then 
$ f $ is not. We  easily see that 
$ \inf_{t \geq 0} F(t)/f(t)=1 $ and that
$E_{0}[\|X\|^{2}]=p(1-ab^{p/2+1})/(1-ab^{p/2}) $.
Because $ (1-ab^{p/2+1})/(1-ab^{p/2}) \leq 2 $ for 
 $ 0< a \leq 1 $, $ 0<b<1 $, the generalized Bayes estimator is minimax
 by \cite{Berger1975}.
\end{example}

The following lemma is stated in a more general form in 3.036 of
\cite{Gradshteyn-Ryzhik}, but it is incorrectly stated with an errata
posted on the book's web page.  Maruyama pointed out this error and he is
acknowledged in the errata for 3.036.  Since a derivation of the formula
is not easily accessible, we provide our own proof.

\begin{lemma}
 \label{int}
For $ \alpha > -1/2 $ and $|a|<1 $,
 \begin{equation*}
  \int_{0}^{\pi}
(1+2a \cos \varphi+a^{2})^{-\alpha}\sin^{2\alpha}\varphi d\varphi=
B(\alpha+1/2,1/2).
 \end{equation*}
\end{lemma}

\begin{proof}
 Let $ g(\varphi)=(1+2a \cos \varphi+a^{2})^{-1}\sin^{2}\varphi $.
Then we have the derivative 
\begin{equation*}
 g'(\varphi)=2 \sin \varphi \frac{(a\cos \varphi +1)(\cos \varphi +a)}{(1+2a \cos \varphi+a^{2})^{2}}.
\end{equation*} 
We see that $ g(\varphi) $ is monotone increasing 
from $ g(0)=0$ to $ g(\arccos (-a))=1 $ and decreasing from
 $g(\arccos (-a))=1 $ to $ g(\pi)=0 $.
Therefore we have
\begin{align}
& \int_{0}^{\pi} (1+2a \cos \varphi+a^{2})^{-\alpha}\sin^{2\alpha}\varphi d\varphi \nonumber \\
& =\left( \int_{0}^{\arccos(-a)}+\int_{\arccos(-a)}^{\pi}\right)
  (1+2a \cos \varphi+a^{2})^{-\alpha}\sin^{2\alpha}\varphi d\varphi \nonumber \\
&=\int_{0}^{\arccos(-a)} (1+2a \cos \varphi+a^{2})^{-\alpha}\sin^{2\alpha}\varphi d\varphi
\nonumber \\
& \quad +\int_{0}^{\arccos(a)} (1-2a \cos \rho+a^{2})^{-\alpha}\sin^{2\alpha}\rho d\rho \nonumber \\
&  = \int_{0}^{1} t^{\alpha}(d\varphi/dt)dt+\int_{0}^{1} s^{\alpha}(d\rho/ds)ds,
 \label{h} 
\end{align}
where $ t=(1+2a \cos \varphi+a^{2})^{-1}\sin^{2}\varphi $ and 
$ s=(1-2a \cos \rho+a^{2})^{-1} \sin^{2}\rho $.
Here $(d\varphi/dt) $ and $ (d\rho/ds) $ are calculated as
\begin{align*}
 d\varphi/dt &=\frac{1}{2tA(t)}\left( 1-\left\{at-A\left(t
 \right)\right\}^2 
\right)^{1/2}\\
 d\rho/ds &=\frac{1}{2sA(s)}\left( 1-\left\{as+A\left(s
 \right)\right\}^2 \right)^{1/2},
\end{align*}
where $ A(t)=(1-t)^{1/2}(1-a^{2}t)^{1/2} $. 
Let
\begin{equation*}
h(t)=\frac{1}{2tA(t)}
\left\{\left( 1-\left\{at-A\left(t \right)\right\}^{2}
\right)^{1/2}+\left( 1-\left\{at+A\left(t \right)\right\}^{2}
\right)^{1/2}\right\}.
\end{equation*}
Then we have
$
h^{2}(t)=\{2tA(t)\}^{-2}\{2-2a^{2}t^{2}-2A(t)^{2}+2B(t)\}
$, 
where
\begin{align*}
 B(t) & =\left( 1-\{at-A(t)\}^{2}-\{at+A(t)\}^{2}
+\{a^{2}t^{2}-A^{2}(t)\}^2\right)^{1/2}\\
& = t(1-a^{2}),
\end{align*}
which implies  $ h(t)=t^{-1/2}(1-t)^{-1/2} $. 
Therefore we get
\begin{equation*}
 \mathrm{the \ right \ hand \ side \ of \
  }\eqref{h}=\int_{0}^{1}t^{\alpha}h(t)dt 
=B(\alpha+1/2,1/2).
\end{equation*}
\end{proof}

\bibliographystyle{ims}
\bibliography{arxiv-bib}

\end{document}